\documentclass[12pt]{amsart}
\usepackage{amssymb}
\usepackage{bm}
\usepackage{graphicx}
\usepackage[centertags]{amsmath}
\usepackage{amsfonts}
\usepackage{amsthm}

\newtheorem{thm}{Theorem}
\newtheorem{cor}[thm]{Corollary}
\newtheorem{lem}[thm]{Lemma}
\newtheorem{prop}[thm]{Proposition}

\newtheorem{defn}[thm]{Definition}
\theoremstyle{definition}

\newcommand{\nn}{\mathbb{N}}
\newcommand{\qq}{\mathbb{Q}}
\newcommand{\ee}{\varepsilon}

\newcommand{\PB}{\mathbf{\Pi}}

\newcommand{\sbs}{\mathrm{SB}}
\newcommand{\refl}{\mathrm{REFL}}
\newcommand{\sd}{\mathrm{SD}}

\newcommand{\sz}{\mathrm{Sz}}
\newcommand{\one}{\mathbf{1}}
\newcommand{\sg}{\sigma}
\newcommand{\spa}{\overline{\mathrm{span}}}

\begin{document}

\title[Strongly bounded classes]{Some strongly bounded classes of Banach spaces}
\author{Pandelis Dodos and Valentin Ferenczi}
\address{National Technical University of Athens, Faculty of Applied
Sciences, Department of Mathematics, Zografou Campus, 157 80,
Athens, Greece} \email{pdodos@math.ntua.gr}
\address{Universit\'e Pierre et Marie Curie - Paris 6,
Equipe d'Analyse Fonctionnelle, Bo\^ite 186, 4, place Jussieu,
75252 Paris Cedex 05, France} \email{ferenczi@ccr.jussieu.fr}

\footnotetext[1]{2000 \textit{AMS classification}: 03E15, 46B03.}
\footnotetext[2]{Research supported by a grant of EPEAEK program
"Pythagoras".}

\maketitle


\begin{abstract}
We show that the classes of separable reflexive Banach spaces and
of spaces with separable dual are strongly bounded. This gives a
new proof of a recent result of E. Odell and Th. Schlumprecht,
asserting that there exists a separable reflexive Banach space
containing isomorphic copies of every separable uniformly convex
Banach spaces.
\end{abstract}


\section{Introduction}
A Banach space $X$ is said to be {\em universal} for a class
$\mathcal{C}$ of Banach spaces when every space in $\mathcal{C}$
embeds isomorphically into $X$. It is {\em complementably}
universal when the embeddings are complemented.

According to the classical Mazur Theorem, the space $C(2^\nn)$ is
universal for separable Banach spaces. By \cite{JS}, there does
not exist a separable Banach space which is complementably
universal for the class of separable Banach spaces. However, A.
Pe\l czy\'nski \cite{P} constructed a space $U$ with a Schauder
basis which is complementably universal for the class of spaces
with a Schauder basis (and even for the class of spaces with the
Bounded Approximation Property). There is also an unconditional
version of $U$, i.e. a space with an unconditional basis which is
complementably universal for the class of spaces with an
unconditional basis.

In 1968, W. Szlenk proved that there does not exist a Banach space
with separable dual which is universal for the class of reflexive
Banach spaces \cite{Sz}. His proof is based on the definition of
the Szlenk index which is a transfinite measure of the
separability of the dual of a space. In 1980, J. Bourgain proved
that any space which is either universal for separable reflexive
spaces, or for all $C(K)$ for $K$ countable compact, must  be
universal for all separable Banach spaces \cite{Bou}. B. Bossard
formalized the use of descriptive set theory initiated by Bourgain
to study classes of separable Banach spaces in \cite{B1,B2}. He
proved that any class of separable Banach spaces which is
analytic, in the Effros-Borel structure of subspaces of
$C(2^\nn)$, and contains all separable reflexive Banach spaces
must contain a universal space. In a recent paper \cite{AD}, S.A.
Argyros and the first named author have connected Bourgain's and
Bossard's approach to universality problems. Among others, they
introduced the following concept.
\begin{defn}
\label{d1} A class $\mathcal{C}$ of separable Banach spaces is
said to be strongly bounded if for every analytic subset $A$ of
$\mathcal{C}$, in the Effros-Borel structure of subspaces of
$C(2^\nn)$, there exists $Y\in\mathcal{C}$ that contains
isomorphic copies of every $X\in A$.
\end{defn}
This notion is central for understanding universality problems in
Banach space theory. In \cite{AD}, it is shown that several
natural classes of separable Banach spaces are strongly bounded.
In particular, the following is proved (see Theorem N in the
introduction of \cite{AD}).
\begin{thm}
\label{t1} The following hold.
\begin{enumerate}
\item[(1)] The class of reflexive spaces with a Schauder basis is
strongly bounded. \item[(2)] The class of spaces with a shrinking
Schauder basis is strongly bounded.
\end{enumerate}
\end{thm}
In this note we remove the assumption of the existence of a basis
in Theorem \ref{t1} and we prove the following.
\begin{thm}
\label{t2} The following hold.
\begin{enumerate}
\item[(1)] The class of separable reflexive spaces is strongly
bounded. \item[(2)] The class of spaces with a separable dual is
strongly bounded.
\end{enumerate}
\end{thm}
Our method is to reduce the proof of Theorem \ref{t2} to Theorem
\ref{t1} by using a uniform version of the Theorem of Zippin
\cite{Z} stating that every Banach space with a separable dual
embeds into a space with a shrinking Schauder basis. To this end,
we are essentially based on the results of B. Bossard in \cite{B1}
and the alternative proof of Zippin's theorem given by N.
Ghoussoub, B. Maurey and W. Schachermayer in \cite{GMS}.

Theorem \ref{t2} answers positively a question of H.P. Rosenthal
from 1979, \cite{Ro}. He asked whether there existed a universal
space with a separable dual for any given class of spaces on which
the Szlenk index is bounded. In particular, we have the following.
\begin{cor}
\label{c1} For every countable ordinal $\xi$, the class of spaces
with Szlenk index less or equal than $\xi$ is Borel. Thus, for
every $\xi<\omega_1$, there exists a Banach space $Y_\xi$ with
separable dual such that for any space $X$ with $\sz(X)\leq\xi$,
$X$ embeds into $Y_\xi$.
\end{cor}
J. Bourgain had asked whether there existed a separable reflexive
Banach space which is universal for separable uniformly convex
spaces. Very recently, E. Odell and Th. Schlumprecht answered this
question by the affirmative \cite{OS}. In our point of view, their
result is an immediate consequence of Theorem \ref{t2} and of the
fact that uniform convexity is a local property, and therefore,
that the class of separable uniformly convex Banach spaces is
Borel.
\begin{cor}[E. Odell, Th. Schlumprecht]
\label{c2} The class $\mathrm{UC}$ of uniformly convex separable
Banach spaces is Borel. Thus, there exists a separable reflexive
Banach space $Y$ that contains isomorphic copies of all uniformly
convex separable Banach spaces.
\end{cor}


\section{Preliminaries}
A topological space is Polish if it is separable and its topology
may be generated by a complete metric. Its Borel subsets are those
belonging to the smallest $\sigma$-algebra containing the open
sets. An analytic subset is the continuous image of a Polish
space, or equivalently, of a Borel subset of a Polish space. A
co-analytic subset is the complement of an analytic subset. If $X$
and $Y$ are Polish spaces, a Borel map $f$ from $X$ into $Y$ is a
map such that $f^{-1}(B)$ is a Borel subset of $X$ for any Borel
subset of $Y$.

If $X$ is a Polish space and $B$ is a co-analytic subset of $X$,
then a map $\phi:B\to\omega_1$ is said to be a co-analytic rank on
$B$ (a $\PB^1_1$-rank in the logical terminology) if there are
relations $\leq_\Sigma, \leq_\Pi$ in $X\times X$ which are
analytic and co-analytic respectively, such that for every $x, y\in
B$ we have
\[ \phi(x)\leq \phi(y) \Leftrightarrow x\leq_\Sigma y \Leftrightarrow
x \leq_\Pi y. \] We refer to \cite{Kechris} for a thorough
presentation of rank theory as well as to \cite{KL} for its
applications. Here we simply state the following properties of
co-analytic ranks which will be needed later on (see
\cite{Kechris}).
\begin{lem}
\label{szl1} Let $X$ be a Polish space, $B$ a co-analytic subset
of $X$ and $\phi:B\to\omega_1$ a co-analytic rank on $B$. Then the
following hold.
\begin{enumerate}
\item[(a)] (Boundedness) For every $A\subseteq B$ analytic we have
\[ \sup\{ \phi(x):x\in A\}<\omega_1. \]
\item[(b)] For every $\xi<\omega_1$, the set $\{ x\in B: \phi(x)\leq\xi\}$
is Borel.
\end{enumerate}
\end{lem}

\subsection{The standard Borel space of separable Banach spaces}
Let $X$ be a Polish space and denote by $F(X)$ the collection of
all closed subsets of $X$. We equip $F(X)$ with the Effros-Borel
$\sg$-algebra. This is the $\sg$-algebra generated by the sets $\{
F\in F(X):F\cap U\neq\varnothing\}$, where $U$ ranges over all
non-empty open subsets of $X$. It is well-known that the
Effros-Borel structure is standard. This means that there exists a
Polish topology $\tau$ on $F(X)$ such that the Borel $\sg$-algebra
of $(F(X),\tau)$ coincides with the Effros-Borel $\sg$-algebra
(see \cite{Kechris}, Theorem 12.6).

Now let $X$ be a separable Banach space and put
\[ \mathrm{Subs}(X)=\{ F\in F(X): F \text{ is a linear subspace of } X\}. \]
Then $\mathrm{Subs}(X)$ is a Borel subset of $F(X)$ (see
\cite{Kechris}, page 79) and so a standard Borel space of its own
right. If $X=C(2^\nn)$, then $\mathrm{Subs}\big(C(2^\nn)\big)$ is
the standard Borel space of all separable Banach spaces and we
denote it simply by $\sbs$. We refer to \cite{AD}, \cite{AGR},
\cite{B0}, \cite{B2} and \cite{Kechris} for more background
material on $\sbs$. We will need the following fact, which is
essentially a consequence of the Kuratowski--Ryll-Nardzewski
selection theorem (see \cite{Kechris}, page 76). There exist two
sequences $d_n:\sbs\to C(2^\nn)$ and $S_n:\sbs\to C(2^\nn)$,
$n\in\nn$, of Borel functions such that for every $X\in\sbs$ we
have $\overline{\{ d_n(X)\}}_n=X$ and
$\overline{\{S_n(X)\}}_n=S_X$.

We denote by $\refl$ and $\sd$ the subsets of $\sbs$ consisting of
all reflexive and all spaces with separable dual respectively.
Both are co-analytic non-Borel (see \cite{B2}). For every
separable space $X$, by $\sz(X)$ we denote the Szlenk index of $X$
(see \cite{Sz}). It is defined as follows. Let $F$ be a
$w^*$-closed subset of $B_{X^*}$. For $\ee>0$, we let $F'_{\ee}$
be the set of $x^*$ in $F$ such that for any $w^*$-neighborhood
$V$ of $x^*$ we have $\mathrm{diam}(V\cap F)>\ee$. Let
$F_{\ee}^{(0)}=B_{X^*}$ and define by transfinite recursion
\[ F_{\ee}^{(\xi)}=(F_{\ee}^{(\zeta)})^{\prime}_{\epsilon}\]
if $\xi=\zeta+1$ is a successor ordinal, and
\[F_{\ee}^{(\xi)}=\bigcap_{\zeta<\xi} F_{\ee}^{(\zeta)}\]
if $\xi$ is a limit ordinal. Then we set
$\sz_{\ee}(X)=\inf\{\xi<\omega_1: F_{\ee}^{(\xi)}=\varnothing\}$
if the set $\{\xi<\omega_1: F_{\ee}^{(\xi)}=\varnothing\}$ is
non-empty, otherwise we set $\sz_{\ee}(X)=\omega_1$. Finally we
let \[ \sz(X)=\sup_{\ee>0} \sz_{\ee}(X).\] It is well-known that
$X\in\sd$ if and only if $\sz(X)<\omega_1$. However, most
important for our purposes is the fact that the Szlenk index is a
co-analytic rank on $\sd$ (see \cite{B2}). Thus Lemma \ref{szl1}
applies to it. For an extensive survey on the Szlenk index we refer
to \cite{L2}.


\section{A uniform version of Zippin's theorem}

The aim of this section is to present the following uniform
version of M. Zippin's theorem \cite{Z} essentially based on the
results of B. Bossard in \cite{B1}.
\begin{prop}
\label{szp1} The following hold.
\begin{enumerate}
\item[(1)] Let $A$ be an analytic subset of $\refl$. Then there
exists an analytic subset $A'$ of $\refl$ such that for all $X\in
A$ there exists $Y\in A'$ with a Schauder basis that contains $X$.
\item[(2)] Let $A$ be an analytic subset of $\sd$. Then there
exists an analytic subset $A'$ of $\sd$ such that every $Y\in A'$
has a shrinking basis and such that for all $X\in A$ there exists
$Y\in A'$ that contains $X$.
\end{enumerate}
\end{prop}
The proof of Proposition \ref{szp1} is modeled after the proof of
Zippin's theorem given by N. Ghoussoub, B. Maurey and W.
Schachermayer in \cite{GMS}.  Part (2) of Proposition \ref{szp1}
is an immediate consequence of the following result of Bossard
(see \cite{B1}, Theorem 3.1) modulo the fact that the Szlenk index
is a co-analytic rank on $\sd$. Note that since having a Schauder
basis is analytic, we may always assume in Proposition \ref{szp1}
that $A'$ is an analytic set of spaces with a Schauder basis.
\begin{thm}[B. Bossard]
\label{szt2} There exists a universal map
$\phi:\omega_1\to\omega_1$ such that for every Banach space $X$
with separable dual and every countable ordinal $\xi$, if
$\sz(X)\leq\xi$, then $X$ embeds into a Banach space $Y$ with a
shrinking basis with satisfies $\sz(Y)\leq\phi(\xi)$.
\end{thm}
To see that Theorem \ref{szt2} implies part (2) of Proposition
\ref{szp1} one argues as follows. Let $A$ be an analytic subset of
$\sd$. By Lemma \ref{szl1}(a), we get
\[ \sup\{ \sz(X):X\in A\}=\xi<\omega_1. \]
Let $[\nn]$ denote the set of all infinite subsets of $\nn$ and let
$(u_n)_n$ denote the basis of the universal space $U$ of Pe\l czy\'nski.
Consider the set
\[ \mathcal{S}=\{ L\in[\nn]: (u_n)_{n\in L} \text{ is shrinking}\}. \]
In \cite{B2}, it is shown that $\mathcal{S}$ is co-analytic and that
the map
\[ \mathcal{S}\ni L\mapsto \sz\big( \overline{\mathrm{span}}\{u_n:n\in L\}\big) \]
is a co-analytic rank on $\mathcal{S}$ (see \cite{B2}, Theorem 5.4).
Therefore, by Lemma \ref{szl1}(b), we get that the set
\[ \mathcal{S}_\xi=\{ L\in\mathcal{S}: \sz\big( \overline{\mathrm{span}}
\{u_n:n\in L\}\big)\leq \phi(\xi)\} \]
is a Borel subset of $\mathcal{S}$. Since the map $[\nn]\ni
L\mapsto
\overline{\mathrm{span}}\{u_n:n\in L\}\in\sbs$ is Borel, it follows that
the set
\[ A'=\{Y\in\sbs: \exists L\in\mathcal{S}_\xi \text{ such that }
\overline{\mathrm{span}}\{u_n:n\in L\}\cong Y\} \]
is an analytic subset of $\sd$ (here $\cong$ denoted as usual the relation
of isomorphism, which is analytic). Theorem \ref{szt2} implies that
$A'$ is as desired.

This simple argument cannot be used in order to derive part (1) of
Proposition \ref{szp1} directly by Theorem \ref{szt2}, as the
Szlenk index is not a co-analytic rank on $\refl$ (see
\cite{B1}, page 68). However, it
does follow from the techniques of \cite{B1} and the method of
\cite{GMS}. We will describe this below. Let $f_0\in C(2^\nn)$ be
a fixed function that separates points in $2^\nn$ and let $\one$
be the constant function equal to $1$. For every $X\in \sbs$ we
let
\[ E(X)=\overline{\mathrm{span}}\{ X\cup f_0\cup \one\}. \]
We have the following easy fact.
\begin{lem}
\label{szl2} The map $\sbs\ni X\mapsto E(X)\in\sbs$
is Borel. In particular, if $A\subseteq \refl$ is analytic, then the
set $A_1=\{ E(X):X\in A\}$ is an analytic subset of $\refl$.
\end{lem}
\begin{proof}
Let $d_n:\sbs\to C(2^\nn)$, $n\in\nn$, be the sequence of Borel
functions such that for all $X\in\sbs$ we have
$\overline{\{d_n(X)\}}_n=X$. Now observe that for every
$U\subseteq C(2^\nn)$ open, we have
\begin{eqnarray*}
E(X)\cap U\neq \varnothing & \Leftrightarrow & \exists p_1,p_2\in\qq \
\exists n_1,...,n_k\in\nn \ \exists q_1,...,q_k\in \qq \\
& & \text{such that } \sum_{i=1}^k q_i d_{n_i}(X)+p_1f_0+p_2\one\in U.
\end{eqnarray*}
Thus the function $E$ is Borel. As for every reflexive space $X$ the space
$E(X)$ is reflexive, the lemma is proved.
\end{proof}
From now on we fix an analytic subset $A$ of $\refl$. Let $A_1$ be
the set obtained by Lemma \ref{szl2} for $A$. Applying Lemma
\ref{szl1}(a) we see that
\[ \sup\{ \sz(Z):Z\in A_1\}=\xi<\omega_1. \]
Denote by $\mathbf{e}=(e_n)_n$ the canonical basis of $\ell_1$. If
$H\in \mathrm{Subs}(\ell_1)$ and $e\in\ell_1$, then $e^H$ will be
the class of $e$ in $\ell_1/H$ and $\mathbf{e}^H=(e_n^H)_n$.
Recall that any separable Banach space is isometric to $\ell_1/H$
for some $H$. By Lemma 3.2 in \cite{B2}, the subset
$\mathcal{Z}_\xi$ of $\mathrm{Subs}(\ell_1)\times \ell_1^\nn\times
C(2^\nn)^\nn\times \sbs$ defined by
\begin{eqnarray*}
\mathcal{Z}_\xi = \{ (H,\mathbf{h},\mathbf{x},X) & : & \sz(X)\leq\xi, \
\spa(\mathbf{x})=X, \ \spa(\mathbf{h})=H \\
& & \mathbf{x}\stackrel{1}{\sim} \mathbf{e}^H, \ \one\in X \text{ and } f_0\in X\}
\end{eqnarray*}
is Borel (as usual, $\mathbf{x}\stackrel{1}{\sim} \mathbf{e}^H$
means that $\mathbf{x}$ is 1-equivalent to $\mathbf{e}^H$). For
every $a\in\mathcal{Z}_\xi$ write
$a=\big(H(a),\mathbf{h}(a),\mathbf{x}(a),X(a)\big)$. Given such
$a$, and following the slicing methods developed in \cite{GMS},
one produces the following.
\begin{enumerate}
\item[(I)] A closed, convex, bounded and symmetric subset $W(a)$
of $C(2^\nn)$ in such a way that the map
\[ \mathcal{Z}_\xi\ni a \mapsto W(a)\in F\big(C(2^\nn)\big) \]
is Borel (see \cite{B2}, Lemma 3.6). Moreover, if $X(a)$ is
reflexive, then the set $W(a)$ is weakly compact. \item[(II)] A
monotone basis $\mathbf{b}(a)\in C(2^\nn)^\nn$ for $C(2^\nn)$ in
such a way again that the map
\[ \mathcal{Z}_\xi\ni a \mapsto \mathbf{b}(a)\in C(2^\nn)^\nn \]
is Borel (see \cite{B2}, Lemma 3.5).
\end{enumerate}
Performing the Davis-Figiel-Johnson-Pe\l czy\'nski interpolation
\cite{DFJP} for the pair $\big( C(2^\nn),W(a)\big)$,  it is shown
in \cite{GMS} that the interpolation space $\Delta(a)$ contains
$X(a)$ and the sequence $\mathbf{b}(a)$ defines a shrinking basis
of $\Delta(a)$. We notice that by (I) above, if $X(a)$ is
reflexive, then classical properties of the interpolation scheme
of \cite{DFJP} imply that the space $\Delta(a)$ is reflexive too.
Denote by $\tilde{\mathbf{b}}(a)$ the sequence $\mathbf{b}(a)$
regarded as a basis of $\Delta(a)$. The crucial fact established
by this procedure is that the subset $R$ of $\mathcal{Z}_\xi\times
C(2^\nn)^\nn\times \sbs$ defined by
\[ R=\{ (a,\mathbf{v},V): \spa(\mathbf{v})=V \text{ and } \mathbf{v}
\stackrel{1}{\sim} \tilde{\mathbf{b}}(a)\}\]
is Borel (see the proof of Lemma 3.7 in \cite{B2}). Notice that if
$(a,\mathbf{v},V)\in R$, then $V$ is isometric to $\Delta(a)$.
Now consider the set $A'$ defined by
\[ V\in A' \Leftrightarrow \exists Y\in A \ \exists a\in \mathcal{Z}_\xi
\ \exists \mathbf{v}\in C(2^\nn)^\nn \ \big[ E(Y)=X(a) \wedge
(a,\mathbf{v},V)\in R\big]. \] As $R$ is Borel and $E$ is a Borel
map, we see that $A'$ is analytic. By the fact that $A\subseteq
\refl$ and property (I) above, we get that every $V\in A'$ is
reflexive, i.e. $A'$ is an analytic subset of $\refl$. Finally, we
notice that for every $Y\in A$ there exists
$a_Y\in\mathcal{Z}_\xi$ such that $X(a_Y)=E(Y)$. Thus
$\Delta(a_Y)\in A'$. As $E(Y)=X(a_Y)$ embeds into $\Delta(a_Y)$,
so does $Y$. It follows that the set $A'$ has all desired
properties and the proof of part (1) of Proposition \ref{szp1} is
completed.


\section{Proof of the main results}

\begin{proof}[Proof of Theorem \ref{t2}]
(1) Let $A$ be an analytic subset of $\refl$. By part (1) of
Proposition \ref{szp1}, there exists an analytic subset $A'$ of
$\refl$ such that for every $Y\in A$ there exists $Z\in A'$ with a
Schauder basis such that $Y$ is contained in $Z$. By Theorem
\ref{t1}(1), the result follows.\\
(2) Let $A$ be an analytic subset of $\sd$. By part (2) of
Proposition \ref{szp1}, there exists an analytic subset $A'$ of
$\sd$ of spaces with a shrinking basis, such that for every $Y\in
A$ there exists $Z\in A'$ such that $Y$ is contained in $Z$. By
Theorem \ref{t1}(2), the result follows.
\end{proof}
\begin{proof}[Proof of Corollary \ref{c1}]
Fix a countable ordinal $\xi$. By Theorem \ref{t2}(2), it is
enough to show that the set $\{X\in\sd: \sz(X)\leq\xi\}$ is Borel.
But this is an immediate consequence of the fact that the Szlenk
index is a co-analytic rank on $\sd$ and Lemma \ref{szl1}(b).
\end{proof}
\begin{proof}[Proof of Corollary \ref{c2}]
By Theorem \ref{t2}(1), it is enough to show that the class
$\mathrm{UC}$ of separable uniformly convex Banach spaces is
Borel. To see this let $S_n:\sbs\to C(2^\nn)$, $n\in\nn$, be the
sequence of Borel functions such that for every $X\in\sbs$ we have
$\overline{\{S_n(X)\}}_n=S_X$. Now observe that
\begin{eqnarray*}
X\in \mathrm{UC} &\Leftrightarrow & \forall k\in\nn \ \exists l\in\nn
\text{ such that } \big[ \forall n,m\in\nn \\
& & \|S_n(X)-S_m(X)\|\geq \frac{1}{k} \Rightarrow
\Big\| \frac{S_n(X)-S_m(X)}{2}\Big\| \leq 1-\frac{1}{l} \big].
\end{eqnarray*}
Invoking the Borelness of the functions $(S_n)_n$ we see that
$\mathrm{UC}$ is Borel.
\end{proof}
It is actually known that a separable space is isomorphic to a
uniformly convex space if and only if its weak*-dentability index
(an ordinal index close to the Szlenk index) is less than $\omega$ \cite{L}.
Our method does not give information on the Szlenk index, or
the weak*-dentability index, of the reflexive space provided by
Corollary C. This space must have weak*-dentability index strictly
greater than $\omega$.

\

We close this section by noting the following unconditional
versions of the above results. Their proofs are easy adaptations
of the methods of \cite{AD}, with the use of the unconditional version of the universal
space of Pe\l czy\'nski instead of the Schauder basis version.

\begin{thm}
\label{unct} (1) The class of reflexive spaces with an
unconditional basis is strong\-ly bounded.\\
(2) The class of spaces with an unconditional basis and not
containing $\ell_1$ is strongly bounded.
\end{thm}

The universal spaces with an unconditional basis and not containing $\ell_1$
obtained in (2) are actually complementably universal for the corresponding
analytic classes of spaces with an unconditional basis and not containing
$\ell_1$.

\begin{cor}
(1) The class $\mathrm{UUC}$ of uniformly convex Banach spaces
with an unconditional basis is analytic. Thus, there exists a
reflexive Banach $Y$ with an unconditional basis that contains
isomorphic copies of all uniformly convex Banach spaces with an
unconditional basis.\\
(2) The class of spaces with an unconditional basis and
non-trivial type is analytic. Therefore there exists a reflexive
space with an unconditional basis which is universal for this
class.
\end{cor}
\begin{proof}
We notice that the class of spaces with an unconditional basis is
analytic. So part (1) is an immediate consequence of Theorem
\ref{unct}(1). For part (2) we recall that the subset of $\sbs$
consisting of all spaces with non-trivial type is analytic (see
\cite{AD}). Observe that a space with non trivial type cannot
contain a copy of $l_1$ or $c_0$, therefore by the classical
theorem of James, it must be reflexive if it has an unconditional
basis. By Theorem \ref{unct}(1), the result follows.
\end{proof}


\end{document}